\documentclass[11pt,letterpaper]{article}

\usepackage[margin=1in]{geometry}

\usepackage{amsmath,amsthm,amssymb}
\usepackage{mathtools}
\usepackage{mathrsfs}
\usepackage{graphicx}
\usepackage{float}
\usepackage{bm}

\usepackage{comment}

\usepackage{hyperref}
\usepackage[nameinlink,capitalise,noabbrev]{cleveref}

\usepackage{tikz}

\theoremstyle{definition}
\newtheorem{definition}{Definition}

\newtheorem{problem}[definition]{Problem}

\theoremstyle{plain}

\newtheorem{theorem}[definition]{Theorem}
\newtheorem{lemma}[definition]{Lemma}
\newtheorem{proposition}[definition]{Proposition}
\newtheorem{claim}[definition]{Claim}

\def \u {\mathbf{u}}
\def \K {\mathcal{K}}

\def \ex {\mathrm{ex}}

\def \ce {\coloneqq}

\renewcommand{\le}{\leqslant}
\renewcommand{\ge}{\geqslant}

\renewcommand{\geq}{\geqslant}

\def \es {\varnothing}

\title{On the maximum number of $r$-cliques in graphs free of \\complete $r$-partite subgraphs}

\author{
J\'ozsef Balogh\thanks{University of Illinois Urbana-Champaign, 1409 W. Green Street, Urbana IL 61801, United States.  Research is supported in part by NSF grants DMS-1764123 and RTG DMS-1937241, FRG DMS-2152488, the Arnold O. Beckman Research Award (UIUC Campus Research Board RB 24012).
Email: \texttt{jobal@illinois.edu}. }
\and
Suyun Jiang\thanks{School of Artificial Intelligence, Jianghan University, Wuhan, Hubei, China, and Extremal Combinatorics and Probability Group (ECOPRO), Institute for Basic Science (IBS), Daejeon, South Korea. Supported by Hubei Provincial Natural Science Foundation of China, National Natural Science Foundation of China (11901246) and China Scholarship Council and IBS-R029-C4.
Email: \texttt{jiang.suyun@163.com}. }
\and
Haoran Luo\thanks{University of Illinois Urbana-Champaign, 1409 W. Green Street, Urbana IL 61801, United States. Research is partially supported by Trjitzinski Fellowship.
Email: \texttt{haoranl8@illinois.edu}. }
}

\date{}

\begin{document}
\maketitle

\begin{abstract}
We estimate the maximum possible number of cliques of size $r$ in an $n$-vertex graph free of a fixed complete $r$-partite graph $K_{s_1, s_2, \ldots, s_r}$. By viewing every $r$-clique as a hyperedge, the upper bound on the Tur\'an number of the complete $r$-partite hypergraphs
gives the upper bound $O\left(n^{r - {1}/{\prod_{i=1}^{r-1}s_i}}\right)$. We improve this  to $o\left(n^{r - {1}/{\prod_{i=1}^{r-1}s_i}}\right)$. The main tool in our proof is the graph removal lemma. We also provide several lower bound constructions.
\end{abstract}

\section{Introduction} \label{sec: intro}
For a graph $H$, the \emph{Tur\'an number} $\ex(n,H)$ is the maximum number of edges an $n$-vertex $H$-free graph can have. To determine the Tur\'an number of various graphs and hypergraphs is a central problem in Extremal Combinatorics and has received considerable attention.

Alon and Shikhelman~\cite{alon2016many} started the systematic study of the following generalized problem related to the Tur\'an number. For graphs $F$ and $H$, the \emph{generalized Tur\'an number} $\ex(n, F, H)$ is the maximum number of copies of $F$ an $n$-vertex $H$-free graph can have. Hence, $\ex(n,H)$ is the case where $F$ is a single edge. Denote by $K_r$ the complete graph on $r$ vertices. Among many other results, Alon and Shikhelman~\cite{alon2016many} proved that if $\chi(H)$, the chromatic number of $H$, is larger than $r$, then $\ex(n, K_r, H) = \Theta(n^r)$.

In this short note, we consider the case where $F = K_r$ and $H = K_{s_1, s_2, \ldots, s_r}$, the complete $r$-partite graph with parts of size $s_1, s_2, \ldots, s_r$. As remarked in~\cite{alon2016many}, the problem of determining the order of $\ex(n, K_3, K_{s,s,s})$ already seems hard; also, to determine $\ex(n, K_3, K_{1,1,2})$ is related to  the triangle removal lemma, where we have $$n^2 \cdot \exp(-O(\sqrt{\log n})) \le \ex(n, K_3, K_{1,1,2}) \le o(n^2)$$ by Behrend's construction~\cite{behrend1946sets} and the triangle removal lemma~\cite{ruzsa1978triple}.

The previous results focus on the case where $r=3$ and $s_1 =1$. Mubayi and Mukherjee~\cite{mubayi2023triangles}  proved that for $ 1 \le a \le b$,
$$
\ex(n, K_3, K_{1,a,b}) = o(n^{3-{1/a}}) \quad \textrm{and} \quad
\Omega(n^2) \le \ex(n,K_3,K_{1,2,2}) \le o(n^{5/2}).
$$
The latter claim was also proved independently by Methuku, Gr\'osz, and Tompkins (unpublished, see~\cite{mubayi2023triangles}) and the first author (see his lecture note~\cite{jozsefnote}).
Methuku (see Theorem 6.1 in~\cite{mubayi2023triangles}) further proved that if $H$ is a graph with $\ex(n,H) = O(n^\alpha)$ for some $\alpha \in (1,2)$, then $\ex(n, K_3, \hat{H}) = o(n^{1+\alpha})$, where $\hat{H}$ is the suspension of $H$, the graph obtained from $H$ by adding a new vertex which is connected to all other vertices.

For every $K_{s_1,s_2,\ldots, s_r}$-free graph $G$, we can construct an $r$-uniform hypergraph on the same vertex set which contains all copies of $K_r$ in $G$ as hyperedges. It is easy to check that this hypergraph is free of $\K^r_{s_1,s_2,\ldots, s_r}$, the complete $r$-uniform $r$-partite hypergraph with parts of size $s_1, s_2, \ldots, s_r$. By the bound on the hypergraph Tur\'an number (see~\cite{erdos1964extremal}),
we have
\begin{equation} \label{equ::exhyp}
    \ex(n, K_r, K_{s_1, s_2, \ldots,s_r}) \le \ex(n, \K_{s_1, s_2, \ldots,s_r}^{(r)}) = O \left(n^{r-{1}/{\prod_{i=1}^{r-1}s_i}} \right).
\end{equation}
Our first result is that we can improve this upper bound as follows.
Independently, Basu, R\"{o}dl, and Zhao~\cite{basu2025number} recently proved similar results.
\begin{theorem} \label{thm: upper bound}
For every positive integer $r \ge 3$ and positive integers $s_1 \le s_2 \le \ldots \le s_r$, we have
$$
\ex(n, K_r, K_{s_1,s_2,\ldots, s_r}) = o\left(n^{r - {1}/{\prod_{i=1}^{r-1}s_i}}\right).
$$
\end{theorem}
\noindent

For the lower bound, we have the following straightforward construction, which generalizes Proposition 3.1 in~\cite{mubayi2023triangles}.
\begin{proposition}\label{thm: lower bound}
For every positive integer $r \ge 3$ and positive integers $s_1 \le s_2 \le \ldots \le s_r$, we have
$$
\ex(n, K_r, K_{s_1,s_2,\ldots, s_r})
\ge
        \left\lfloor \frac{n}{2} \right\rfloor \cdot
        \ex \left(\left\lceil \frac{n}{2} \right\rceil, K_{r-1}, K_{s_1,s_2,\ldots, s_{r-1}} \right).
    $$
\end{proposition}
\noindent
For example, if $\ex(n, K_{s_1,s_2})= \Omega(n^{2 - {1}/{s_1}})$, then by iterating \cref{thm: lower bound}, we have
$$
    \ex(n,K_r,K_{s_1,s_2,\ldots,s_r}) = \Omega \left(n^{r- {1}/{s_1}} \right).
$$
By known bounds on the Tur\'an numbers~\cite{erdos1966problem, brown1966graphs, kollar1996norm, alon1999norm} of bipartite graphs, we have
$$
\ex(n, K_3, K_{2,2,2}) = \Omega(n^{5/2}), \quad
\ex(n, K_3, K_{3,3,3}) = \Omega(n^{8/3}),
$$
and when $s_3 \ge s_2 \ge (s_1-1)! + 1$, we have
\begin{equation} \label{equ: prop2K3}
\ex(n, K_3, K_{s_1, s_2, s_3}) = \Omega\left(n^{3- {1}/{s_1}} \right).
\end{equation}

It turns out that when $s_r$ is much larger than the other $s_i$'s, we have the following better lower bound.
\begin{proposition} \label{thm: better lower bound}
    For positive integers $s_1 \le s_2 \le s_3$, where $s_3 > (s_1+s_2-1)!$, we have
    $$
        \ex(n,K_3,K_{s_1,s_2,s_3}) \ge \ex(n,K_3,K_{s_1+s_2,s_3})
        = \Omega \left(n^{3-{3}/{(s_1+s_2)}}\right).
    $$
In general,
    for positive integers $s_1 \le s_2 \le \ldots\le s_r$, where $s_r > \left(\sum_{i=1}^{r-1}s_i-1\right)!$ and $ \sum_{i=1}^{r-1}s_i\geq 2r-2$, we have
    $$
        \ex(n,K_r,K_{s_1,s_2,\ldots, s_r}) \ge \ex(n,K_r,K_{\sum_{i=1}^{r-1}s_i,s_r})
        = \Omega \left(n^{r-{\binom{r}{2}}/{\sum_{i=1}^{r-1} s_i}} \right).
    $$
\end{proposition}
\noindent
For example, when $s_3 \gg s_2 \ge 2s_1$, \cref{thm: better lower bound} gives a better bound than the one in~\eqref{equ: prop2K3}.

\section{Proofs}

\subsection{Upper bound}
The previous methods in~\cite{mubayi2023triangles} seem not to work for the cases where $r > 3$ or $s_1 > 1$ without new ideas, as in these cases it is not enough to only consider the degree or neighborhood of some vertices. Here, we take an alternative approach by considering the number of copies of a sequence of graphs, where every graph is obtained by taking the blow-up of the previous one and the last one is $K_{s_1,s_2,\ldots, s_r}$, the graph we forbid.

One of our main tools for \cref{thm: upper bound} is the graph removal lemma, see~\cite{ruzsa1978triple, erdos1986asymptotic} and also~\cite{alon1994algorithmic, furedi1995extremal}.

\begin{lemma} \label{lem: graph removal lemma}
For a graph $H$, every graph on $n$ vertices with $o(n^{|V(H)|})$ copies of $H$ can be made $H$-free by removing $o(n^2)$ edges.
\end{lemma}

We make use of the following lemma about the number of the blow-ups of a given graph. We find it more convenient to talk about the ordered copies of a graph.
Note that if the number of copies and ordered copies of $H$ in $G$ are $N$ and $M$, respectively, then trivially we have
$$
N \le M \le |V(H)|!\cdot N.
$$
For integers $a,b \ge 0$, we let $(a)_b \ce \prod_{i=0}^{b-1} (a-i)$.
\begin{lemma} \label{lem::numBlowup}
    Let $F$ be a $k$-vertex graph, $v$ be a vertex in $V(F)$ and $H$ be the graph obtained from $F$ by adding $t-1$ new copies of $v$. For every $n$-vertex graph $G$, if the number of ordered copies of $F$ in $G$ is $x \cdot n^{k-1}$, where $x \ge 4t$, then the number of ordered copies of $H$ in $G$ is at least
    $\frac{1}{2^t} \cdot x^t \cdot n^{k-1}.$
\end{lemma}
\begin{proof}
Without loss of generality, we assume that $V(F) = \{v_1,\ldots,v_k\}$ and $v = v_k$. Let $V^{k-1}(G)$ be the set of (ordered) sequences of $k-1$ distinct vertices in $G$. For $\u = (u_1, \ldots, u_{k-1}) \in V^{k-1}(G)$, let $M_{\u}$ be the number of ordered copies of $F$ in $G$ containing $u_1, \ldots, u_{k-1}$, in which $u_i$ acts as $v_i$ for $1\le i \le k-1$. Let $A$ be the collection of $\u \in V^{k-1}(G)$ with $ M_{\u} \ge t$.
Now, the number of ordered copies of $H$ in $G$ is
    \begin{align*}
        \sum_{\u \in V^{k-1}(G)}  \left( M_{\u}  \right)_t
         \ge &
        \sum_{\u \in A} \left( M_{\u} - t + 1  \right)^t
         \ge
        \frac{1}{|A|^{t-1}} \left( \sum_{\u \in A}
        (M_{\u} - t + 1)
        \right)^t \\
         \ge &
         \frac{1}{n^{(k-1)(t-1)}} \left( \sum_{\u \in A}
        (M_{\u} - t + 1)
        \right)^t
         \ge
        \frac{1}{n^{(k-1)(t-1)}} \left( \sum_{\u \in V^{k-1}(G)}
        M_{\u} - 2t n^{k-1}
        \right)^t \\
        \ge &
        \frac{1}{n^{(k-1)(t-1)}} \cdot \left(\frac{x}{2} \cdot n^{k-1}\right)^t
        = \frac{1}{2^t} \cdot x^t \cdot n ^ {k-1}
        . \qedhere
    \end{align*}
\end{proof}

\begin{proof}[Proof of \cref{thm: upper bound}]
Let $G$ be an $n$-vertex $K_{s_1,s_2,\ldots, s_r}$-free graph with $N \ce \ex(n,K_r, K_{s_1,s_2,\ldots, s_r})$ copies of $K_r$.
For $i = 1,\ldots, r$, let $M_i$ be the number of ordered copies of the graph $H_i \ce K_{s_1,\ldots, s_i, 1,\ldots, 1}$, the complete $r$-partite graph with $i$ parts of size $s_1, \ldots, s_i$ respectively and $r-i$ parts of size $1$, in $G$.
Note that $H_r$ is just $K_{s_1,s_2,\ldots, s_r}$ and $H_{i}$ can be obtained from $H_{i-1}$ by adding $s_i -1$ new copies of the vertex in a part of size $1$.
We first have the following claim on the upper bounds on $M_i$. For simplicity, we let $\prod_{j=r}^{r-1} s_j = 1$.

\begin{claim} \label{cla::Mi+1Mi}
For $i = 1, \ldots, r-1$, we have
$$
    M_i = O \left( n^{|V(H_i)|- 1/\prod_{j = i+1}^{r-1} s_j} \right).
$$
\end{claim}
\begin{proof}
For each $i= 1, \ldots, r-2$, we apply \cref{lem::numBlowup} with $(F, H, k,t)=(H_{i}, H_{i+1}, |V(H_{i})|, s_{i+1})$. By assumption, we have $M_r = 0$ and hence
$
    M_{r-1} = O(n^{|V(H_{r-1})|-1}).
$
For every $i < r-1$, if
$$
M_{i+1}
= O \left( n^{|V(H_{i+1})|- 1/\prod_{j = i+2}^{r-1} s_j}\right)
= O \left( n^{|V(H_i)|- 1} \cdot n^{s_{i+1} - 1/\prod_{j = i+2}^{r-1} s_j} \right),
$$
then we have
\[
M_i
= O\left( n^{|V(H_i)| - 1} \cdot n^{1 - 1/\prod_{j=i+1}^{r-1}s_j} \right)
= O \left( n^{|V(H_i)|- 1/\prod_{j = i+1}^{r-1} s_j} \right)
.
\]
The claim is proved.
\end{proof}

By \cref{cla::Mi+1Mi}, we have $M_1 = O(n^{s_1 + r-1 -{1}/{\prod_{j=2}^{r-1}s_j}})$.
By~\eqref{equ::exhyp}, we have $N = o(n^r)$, so by~\cref{lem: graph removal lemma}, there is $S \subseteq E(G)$ with $|S| = o(n^2)$ such that removing $S$ makes $G$ free of $K_r$. Let $F_1, F_2, \ldots, F_a$ be the copies of $K_{r-1}$ containing any edge in $S$, where $a \le |S| \cdot n^{r-1 -2} = o(n^{r-1})$.
For $i=1 ,\ldots, a$, let $f_i$ be the number of copies of $K_r$ that $F_i$ is in. Let $B$ be $\{i: 1 \le i \le a, \, f_i \ge s_1\}$.
Note that every copy of $K_r$ in $G$ contains at least one edge in $S$, and hence $\sum_{i=1}^a f_i \ge N$. We may also assume that $N \ge 4s_1a$, since otherwise $N = o(n^{r-1})$ and we are done.
Now, we have
\begin{align*}
O\left(n^{s_1+r-1 -{1}/{\prod_{j=2}^{r-1}s_j}}\right) &= M_1 \ge \sum_{i=1}^a (f_i)_{s_1} \ge \sum_{i \in B} (f_i - s_1 + 1)^{s_1}
\\
&\ge \frac{1}{|B|^{s_1-1}} \left(\sum_{i \in B} (f_i- s_1 + 1)\right)^{s_1}
\ge
\frac{1}{a^{s_1-1}} \left( \sum_{i=1}^a f_i - 2s_1a \right)^{s_1}
\\
&\ge
\frac{1}{a^{s_1-1}} \left( \frac{\sum_{i=1}^a f_i}{2} \right)^{s_1}
\ge \frac{1}{2^{s_1} a^{s_1-1}} N^{s_1},
\end{align*}
and hence
\begin{align*}
    N^{s_1} \le 2^{s_1}a^{s_1-1} \cdot O\left(n^{s_1+r-1 -{1}/{\prod_{j=2}^{r-1}s_j}}\right)
    &= o\left(n^{(r-1)(s_1-1)}\right) \cdot O\left(n^{s_1 + r-1 -{1}/{\prod_{j=2}^{r-1}s_j}}\right)  \\
    &= o\left(n^{rs_1 -{1}/{\prod_{j=2}^{r-1}s_j}} \right).
\end{align*}
Therefore, we have
\begin{equation*}
N = o \left(n^{r - {1}/{\prod_{j=1}^{r-1}s_j}}\right). \qedhere
\end{equation*}
\end{proof}

\subsection{Lower bounds}

\begin{proof}[Proof of \cref{thm: lower bound}]
Let $G$ be a graph where $V(G)$ is partitioned into $(S_1, S_2)$ such that $G[S_1]$
is an $\left\lceil \frac{n}{2} \right\rceil$-vertex
$K_{s_1,s_2,\ldots, s_{r-1}}$-free graph with $\ex \left(\left\lceil \frac{n}{2} \right\rceil , K_{r-1}, K_{s_1,s_2,\ldots, s_{r-1}}\right)$ copies of $K_{r-1}$, $G[S_1, S_2]$ is complete,
and $S_2$ is an independent set.
Then, every copy of $K_{r-1}$ in $S_1$ and every vertex in $S_2$ form a copy of $K_r$,
so $G$ has at least $\left\lfloor \frac{n}{2} \right\rfloor \cdot
        \ex \left(\left\lceil \frac{n}{2} \right\rceil, K_{r-1}, K_{s_1,s_2,\ldots, s_{r-1}}\right)$ copies of $K_r$.

Now we need to prove that $G$ is $K_{s_1,s_2,\ldots, s_{r}}$-free. Suppose that $H$ is a subgraph of $G$ such that $H\cong K_{s_1,s_2,\ldots, s_{r}}$. Since $G[S_1]$ is $K_{s_1,s_2,\ldots, s_{r-1}}$-free and $K_{s_1,s_2,\ldots, s_{r-1}} \subset K_{s_1,s_2,\ldots, s_{r}}$, we have $V(H)\cap S_2 \neq \es$. Assume that $v \in V(H) \cap S_2$. Then, by our assumption that $s_1 \le s_2 \le \ldots \le s_r$ and $S_2$ is an independent set, we have $G[N(v)] \subseteq G[S_1]$ contains a copy of $K_{s_1,s_2,\ldots, s_{r-1}}$, a contradiction.
\end{proof}

To prove \cref{thm: better lower bound} we need the following lemmas.
\begin{lemma}[Lemma 4.4 in \cite{alon2016many}]\label{lem: triangle in bipartite graphs}
    For integers $s\geq 2$ and $t\geq (s-1)!+1$, we have $$\ex(n, K_3, K_{s,t})=\Theta\left(n^{3-{3}/{s}}\right).$$
\end{lemma}
\begin{lemma}[Theorem 1.3 in \cite{alon2016many}]\label{lem: Km in bipartite graphs}
    For integers $m$ and $t\geq s$ satisfying $s\geq 2m-2$ and $t\geq (s-1)!+1$, we have $$\ex(n,K_m,K_{s,t})=\Theta\left(n^{m-{\binom{m}{2}}/{s}}\right).$$
\end{lemma}

\begin{proof}[Proof of \cref{thm: better lower bound}]
For the $K_3$ case,
    note that $K_{s_1+s_2,s_3} \subseteq K_{s_1,s_2,s_3}$, so
    the first inequality follows from the definition.
    By \cref{lem: triangle in bipartite graphs}, the second inequality holds.

The general case follows from the observation $K_{\sum_{i=1}^{r-1}s_i,s_r} \subseteq K_{s_1,s_2,\ldots, s_r}$ and \cref{lem: Km in bipartite graphs}.
\end{proof}

\section{Concluding remarks}
The bounds in \cref{thm: upper bound} can also be applied to get some bounds for $\ex(n,K_r, H)$ where $\chi(H) < r$. For example, noting that in a $K_{1,2,2}$-free graph, every copy of $K_3$ can be in at most one copy of $K_4$, we have $\ex(n, K_4, K_{1,2,2}) \le \ex(n, K_3, K_{1,2,2}) = o(n^{5/2})$. In general, we have
$$
\ex(n, K_{r+1}, K_{1, \ldots, 1 ,2 ,2}) \le \ex(n, K_r, K_{1, \ldots, 1, 2, 2}) \le o(n^{r- 1/2})
.
$$

As mentioned in \cref{sec: intro}, the problem of determining the correct order of $\ex(n,K_r, K_{s_1,s_2,\ldots, s_r})$ seems to be hard in general. In particular, if the upper bound given in \cref{thm: upper bound} is closer to the truth, which we believe, then a better construction than the one in \cref{thm: lower bound} is needed. However, as $\ex(n, K_r, K_{s_1,s_2,\ldots, s_r}) \le \ex(n, \K^{(r)}_{s_1,s_2,\ldots, s_r})$, this will require us to improve the lower bound for the Tur\'an numbers of complete $r$-partite $r$-uniform hypergraphs, which is a notoriously difficult problem.
There are several cases where the orders of $\ex(n, \K^{(r)}_{s_1,s_2,\ldots, s_r})$ are determined (see~\cite{mubayi2002some} and~\cite{ma2018some}), so it may be easier to improve the lower bounds in these cases. We conclude this note with the following problem, which is also emphasized in~\cite{mubayi2023triangles}.

\begin{problem}
    Is it true that
    $$
        \ex(n, K_3, K_{1,2,2}) \ge n^{5/2 - o(1)} \textrm{ ?}
    $$
\end{problem}

\section*{Acknowledgment}
This work was partially done when the first and third authors visited the Institute for Basic Science (IBS). We are very grateful for the kind hospitality of IBS and other visitors who participated in
fruitful discussions at the beginning of this project.
We would also like to thank the anonymous referees for their helpful comments.

\end{document}